\documentclass[12pt,twoside]{article}
\usepackage{graphicx}

\newtheorem{theorem}{\bf Theorem}[section]

\newtheorem{lemma}[theorem]{\bf Lemma}
\newtheorem{corollary}[theorem]{\bf Corollary}

\newcommand{\D}{\Delta}
\newcommand{\de}{\delta}
\newcommand{\To}{\longrightarrow}

\input{amssym}

\date{September 10, 2009}
\begin{document}
\title{{\large On edge-group choosability of graphs}}
\author{\small A. Khamseh$^{\textrm{a}}$, G.R. Omidi$^{\textrm{a},\textrm{b},1}$\\
\small  $^{\textrm{a}}$Department of Mathematical Sciences, Isfahan
University
of Technology,\\ \small Isfahan, 84156-83111, Iran\\
\small  $^{\textrm{b}}$School of Mathematics, Institute for Research
in Fundamental Sciences (IPM),\\
\small  P.O.Box:19395-5746, Tehran,
Iran\\
\small \texttt{E-mails:khamseh@math.iut.ac.ir, romidi@cc.iut.ac.ir}}
\date {}

\maketitle\footnotetext[1] {\tt This research was in part supported
by a grant from IPM (No.88050012)} \vspace*{-0.5cm}

\begin{abstract}

 In this paper, we study the concept of {\it{edge-group
choosability}} of graphs. We say that $G$ is edge-$k$-group
choosable if its line graph is $k$-group choosable. An edge-group
choosability version of Vizing conjecture is given. The evidence of
our claim are graphs with maximum degree less than $4$, planar
graphs with maximum degree at least $11$, planar graphs without
small cycles, outerplanar graphs and near-outerplanar graphs.

\end{abstract}
\section{\normalsize\bf{Introduction}}
We consider only simple graphs in this paper unless otherwise
stated. For a graph $G$, we denote its vertex set, edge set, minimum
degree, and maximum degree by $V(G)$, $E(G)$, $\de(G)$, and $\D(G)$,
respectively. A {\it plane graph} is a particular drawing of a
planar graph in the Euclidean plane. We denote the set of faces of a
plane graph $G$ by $F(G)$. For a plane graph $G$ and $f \in F(G)$,
we write $f = [u_1u_2\ldots u_n]$ if $u_1, u_2, \ldots, u_n$ are the
vertices on the boundary walk of $f$ enumerated clockwise. The
{\it{degree of a face}} is the number of edge-steps in the boundary
walk. Let $d_G(x)$, or simply $d(x)$, denote the degree of a vertex
(or face) $x$ in $G$. A vertex (or face) of degree $k$ is called a
{\it{$k$-vertex}} (or {\it{$k$-face}}). A $3$-face $f = [u_1u_2u_3]$
is called a {\it{$(i,j,k)$-face}} if
$(d(u_1),d(u_2),d(u_3))=(i,j,k)$. For $v\in V(G)$, $N_G(v)$ is the
set of all vertices of $G$ that are adjacent to $v$ in $G$. A
{\it{$2$-alternating cycle}} in a graph $G$ is a cycle of even
length in which alternate vertices have degree $2$.

A {\it{$k$-coloring}} of a graph $G$ is a mapping $\phi$ from $V(G)$
to the set of colors $\{1, 2, \ldots, k\}$ such that $\phi(x)\neq
\phi(y)$ for every edge $xy$. A graph $G$ is {\it{$k$-colorable}} if
it has a $k$-coloring. The {\it{chromatic number}} $\chi(G)$ is the
smallest integer $k$ such that $G$ is $k$-colorable. A mapping $L$
is said to be a {\it{list assignment}} for $G$ if it supplies a list
$L(v)$ of possible colors to each vertex $v$. A {\it{$k$-list
assignment}} of $G$ is a list assignment $L$ with $|L(v)|=k$ for
each vertex $v\in V(G)$. If $G$ has some $k$-coloring $\phi$ such
that $\phi(v) \in L(v)$ for each vertex $v$, then $G$ is
{\it{$L$-colorable}} or $\phi$ is an {\it{$L$-coloring}} of $G$. We
say that $G$ is {\it{$k$-choosable}} if it is $L$-colorable for
every $k$-list assignment $L$. The {\it{choice number}} or {\it{list
chromatic number}} $\chi_l(G)$ is the smallest $k$ such that $G$ is
$k$-choosable. To distinguish the objects by different notions we
denote the line graph of a graph $G$ by $\ell(G)$. By considering
colorings for $E(G)$, we can define analogous notions such as
{\it{edge-$k$-colorability}}, {\it{edge-$k$-choosability}}, the
{\it{chromatic index}} $\chi'(G)$, the {\it{choice index}}
$\chi'_l(G)$, etc. Clearly, we have $\chi'(G)=\chi(\ell(G))$ and
$\chi'_l(G)=\chi_l(\ell(G))$. The notion of list coloring of graphs
has been introduced by Erd\H{o}s, Rubin, and Taylor \cite{2ofmain}
and Vizing \cite{13ofmain}. The following conjecture, which first
appeared in \cite{3of12}, is
well-known as the List Edge Coloring Conjecture.\\[1pt]

{\bf {Conjecture 1.}} If $G$ is a multi-graph, then
$\chi'_l(G)=\chi'(G)$.\\[1pt]

Although Conjecture $1$ has been proved for a few special cases such
as bipartite multigraphs \cite{galvin}, complete graphs of odd order
\cite{5of8}, multicircuits \cite{16of8}, graphs with $\D(G) \geq 12$
that can be embedded in a surface of non-negative characteristic
\cite{14}, and outerplanar graphs \cite{1}, it is regarded as very
difficult. Vizing proposed the following
weaker conjecture (see \cite{10of6}).\\[1pt]

{\bf {Conjecture 2.}} Every graph $G$ is
edge-$(\D(G)+1)$-choosable. \\[1pt]

Assume $A$ is an Abelian group and $F(G, A)$ denotes the set of all
functions $f : E(G)\To A$. Consider an arbitrary orientation of $G$.
The graph $G$ is {\it{$A$-colorable}} if for every $f \in F(G, A)$,
there is a vertex coloring $c : V (G) \To A$ such that $c(x) -c(y)
\neq f(xy)$ for each directed edge from $x$ to $y$. The {\it{group
chromatic number}} of $G$, $\chi_g(G)$, is the minimum $k$ such that
$G$ is $A$-colorable for any Abelian group $A$ of order at least
$k$. The notion of group coloring of graphs was first introduced by
Jaeger et al. \cite{3ofmain}.

The concept of {\it {group choosability}} is introduced by Kr\'{a}l
and Nejedl\'{y} \cite{kral}. Let $A$ be an Abelian group of order at
least $k$ and $L : V (G) \To 2^A$ be a list assignment of $G$. For
$f \in F(G, A)$, an {\it{$(A, L, f)$-coloring}} under an orientation
$D$ of $G$ is an $L$-coloring $c : V (G) \To A$ such that $c(x)-c(y)
\neq f(xy)$ for every edge $e = xy$, $e$ is directed from $x$ to
$y$. If for each $f \in F(G, A)$ there exists an $(A, L,
f)$-coloring for $G$, then we say that $G$ is {\it{$(A,
L)$-colorable}}. The graph $G$ is $k$-group choosable if $G$ is $(A,
L)$-colorable for each Abelian group $A$ of order at least $k$ and
any $k$-list assignment $L : V (G) \To {A\choose k}$. The minimum
$k$ for which $G$ is $k$-group choosable is called the {\it{group
choice number}} of $G$ and is denoted by $\chi_{gl}(G)$. It is clear
that the concept of group choosability is independent of the
orientation on $G$.

Graph $G$ is called {\it{edge-$k$-group choosable}} if its line
graph is $k$-group choosable. The {\it{group-choice index}} of $G$,
$\chi'_{gl}(G)$, is the smallest $k$ such that $G$ is edge-$k$-group
choosable, i.e. $\chi'_{gl}(G)=\chi_{gl}(\ell(G))$. It is easily
seen that an even cycle is not edge-$2$-group choosable. This
example shows that $\chi'_{gl}(G)$ is not generally equal to
$\chi'(G)$. But we can extend the Vizing Conjecture
as follows.\\[1pt]

{\bf {Conjecture 3.}} If $G$ is a multi-graph, then
$\chi'_{gl}(G)\leq\D(G)+1$.\\[1pt]

Since $\D(G)\leq \chi'(G)\leq \D(G)+1$, as a sufficient condition,
we have the following weaker conjecture.\\[1pt]

{\bf {Conjecture 4.}} If $G$ is a multi-graph, then
$\chi'_{gl}(G)\leq\chi'(G)+1$.\\[1pt]

In this paper, we prove that Conjecture $3$ (and consequently
Conjecture $4$) holds for certain classes of graphs. We shall see
that Conjecture $3$ is true for graphs with $\D(G)\leq 3$, planar
graphs with $\D(G)\geq 11$, planar graphs without small cycles,
outerplanar graphs, and near-outerplanar graphs.
\section{\normalsize\bf{Edge-group choosability of graphs with bounded
degree}}\label{secdelta}

In this section, we aim to prove Conjecture $3$ holds for a graph
$G$ satisfying either $\D(G)\leq 3$ or planar with $\D(G)\geq 11$.

\begin{lemma}\label{lem1} Let $l$ be a natural number, $v$ be a vertex of degree
at most $2$ of $G$ and $e$ be an edge adjacent to $v$. If
$\chi'_{gl}(G-e)\leq \Delta(G)+l$, then $\chi'_{gl}(G)\leq
\Delta(G)+l$.
\end{lemma}
{\it Proof.} Let $\D=\D(G)$, $D$ be an orientation  of $\ell(G)$,
$A$ be an Abelian group of order at least $\D+l$, $L:V(\ell(G))\To
{A\choose \D+l}$ be a $(\D+l)$-list assignment and $f\in
F(\ell(G),A)$. Suppose that $G'=G-e$. Then $\ell(G')=\ell(G)-e$ and
since $\chi'_{gl}(G')\leq \Delta+l$, there exists an $(A, L,
f)$-coloring $c:V(\ell(G'))\To A$. For each $e'\in N_{\ell(G)}(e)$
we can consider, without loss of generality, $ee'$ to be directed
from $e$ to $e'$. Then, since $|L(e)|=\D+l$ and $d_{\ell(G)}(e)\leq
\D$, $|L(e)-\{c(e')+f(ee'):e'\in N_{\ell(G)}(e)\}|\geq 1$ and so
there is now at least one color available for $e$. Thus we can color
all edges of $G$. This completes the proof of lemma.$
\hfill \dashv $ \\

An argument similar to the proof of Lemma \ref{lem1} gives the
following lemma.
\begin{lemma}\label{mohem} Let $G$ be a graph with
$\chi'_{gl}(G-e)< \chi'_{gl}(G)$ for each $e\in E(G)$. Then
$\delta(\ell(G))\geq \chi'_{gl}(G)-1$.
\end{lemma}

\begin{lemma}\label{brookscycle}\cite{main} Let $P_n$ and $C_n$ denote a path
and a cycle of length $n$, respectively. Then
\begin{itemize}\item[$(1)$]  $\chi_{gl}(P_n)=2$ and
$\chi_{gl}(C_n)=3$,
\item[$(2)$] For any connected simple graph $G$, we have
$\chi_{gl}(G)\leq \D(G)+1$, with equality holds iff $G$ is either
a cycle or a complete graph.
\end{itemize}
\end{lemma}

Immediately from Lemma \ref{brookscycle}, we have the following
corollary.

\begin{corollary}\label{edgecycle}
$\chi'_{gl}(P_n)=\D(P_n)=2$ and $\chi'_{gl}(C_n)=\D(C_n)+1=3$.
\end{corollary}

\begin{theorem}\label{dless4}
Let $G$ be a graph with maximum degree $\D(G)$. If $\D(G)\leq 3$,
then $\chi'_{gl}(G)\leq \Delta(G)+1$ and if $\D(G)=4$, then
$\chi'_{gl}(G)\leq6$.
\end{theorem}
{\it Proof.} It is clear that we can assume $G$ is connected. If
$\D(G)=1$, then $G=P_2$ and this theorem trivially holds. If
$\D(G)=2$, then $G=P_n$ or $G=C_n$ and the assertion holds by
Corollary \ref{edgecycle}. It is clear that $\D(\ell(G))\leq 4$ if
$\D(G)\leq 3$ and $\D(\ell(G))\leq 6$ if $\D(G)\leq 4$. The proof is
completed by Lemma \ref{brookscycle}. $
\hfill \dashv $ \\

The remainder of this paper consists of an investigation of
Conjecture $3$ for certain classes of planar graphs. Our first
result shows that this conjecture is true for planar graphs with
maximum degree at least $11$.

\begin{lemma}\label{naserasr}\cite{naserasr}
For every planar graph $G$ with minimum degree at least $3$ there is
an edge $e=uv$ with $d(u)+d(v)\leq 13$.
\end{lemma}

\begin{lemma}\label{lem2}
If $G$ is a planar graph with maximum degree $\D$, then
$$\chi'_{gl}(G)\leq \max\{\D+1,12\}.$$
\end{lemma}
{\it Proof.} Let $G$ be a minimal counterexample to Lemma
\ref{lem2}. By Lemmas \ref{lem1} and \ref{naserasr}, there exists
$e\in V(\ell(G))$ with $d_{\ell(G)}(e)\leq 11$, which is a contradiction by Lemma \ref{mohem}.$\hfill \dashv $ \\

The truth of Conjecture $3$ for planar graphs with maximum degree at
least $11$ immediately follows from Lemma \ref{lem2}. In fact we
have the following.
\begin{corollary}
Let $G$ be a planar graph with maximum degree $\D$.
\begin{itemize}
\item[$(1)$] If $\D\geq 11$, then $\chi'_{gl}(G)\leq \D+1$,
\item[$(2)$] If $\D\geq 10$, then $\chi'_{gl}(G)\leq \D+2$.
\end{itemize}
\end{corollary}
\section{\normalsize\bf{Edge-group choosability of planar graphs without small cycles}}\label{seccycles}

\begin{figure}[h]
  \begin{center}
  \includegraphics[width=3cm]{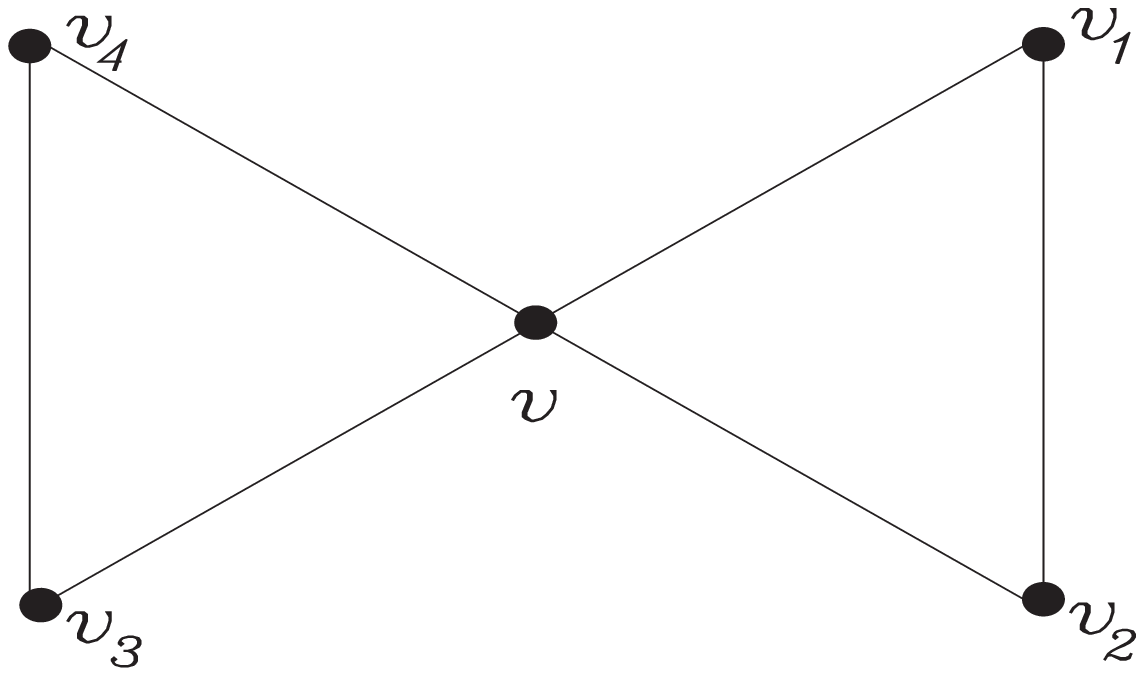}\hspace*{2cm}
  \includegraphics[width=4.5cm]{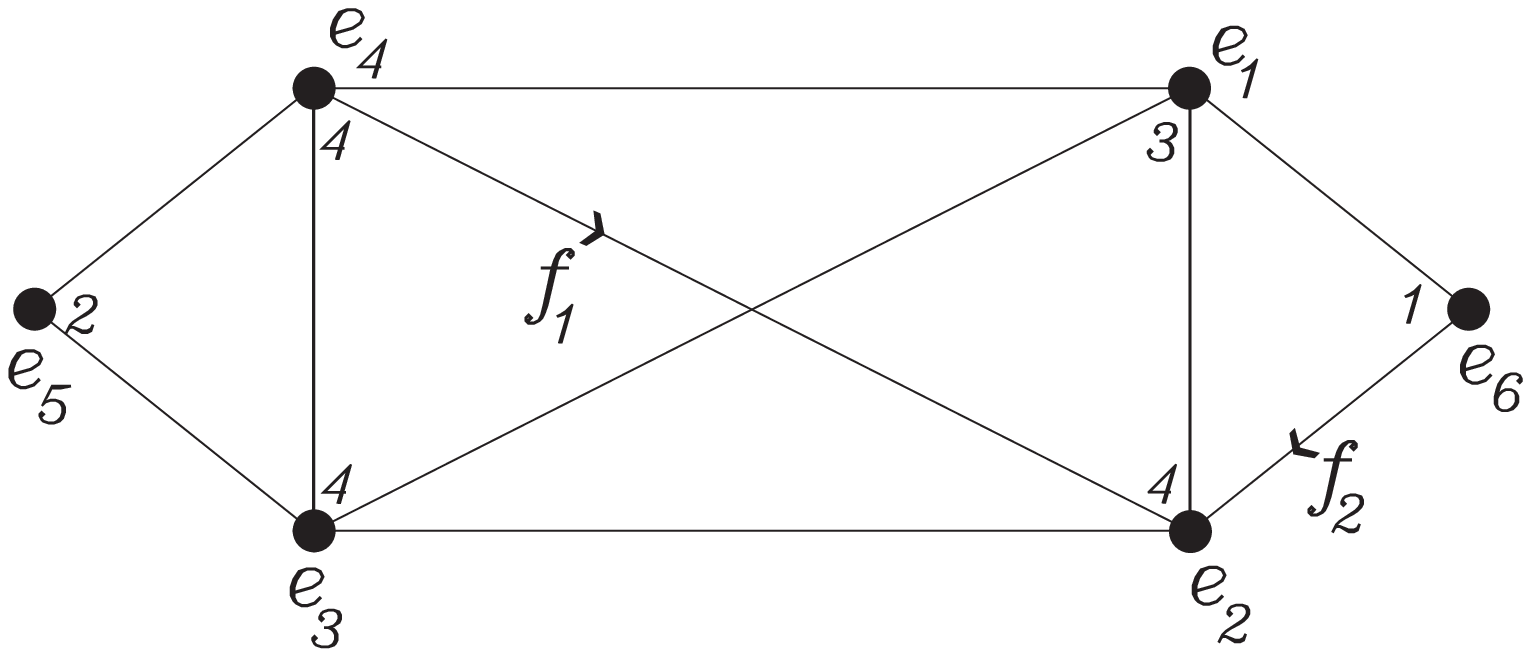}\\
  \caption{The graph $G_1$ and its line graph.}\label{G1}
  \end{center}
\end{figure}

In this section, we show that Conjecture $3$ holds  for planar
graphs that contain no certain configurations. For planar graphs
without $4$-cycles we need the following structural lemma.
\begin{lemma}\label{4-cycles1}\cite{6}
If a connected plane graph $G$ with $\de(G)\geq 3$ has no
$4$-cycles, then $G$ contains one of the following configurations.
\begin{itemize}
\item[$(1)$] An edge $xy$ with $d(x)+d(y)\leq 7$, \item[$(2)$] The
subgraph $G_1$ consisting of a $4$-vertex $v$ that is incident to
two non-adjacent $3$-faces $f_1=[vv_1v_2]$ and $f_2=[vv_3v_4]$, such
that $d(v_i)\geq 4$ for $i=1,2,3,4$, and in $\{v_1,v_2,v_3,v_4\}$
there is at most one vertex with degree greater than $4$ (see Figure
\ref{G1}).
\end{itemize}
\end{lemma}

For the graph $G_1$ of Lemma \ref{4-cycles1}, set $e_i=vv_i$, for
$i=1,2,3,4$, $e_5=v_3v_4$ and $e_6=v_1v_2$. Now we have the
following result.
\begin{lemma}\label{4-cycles2}
Suppose that $A$ is an Abelian group with $|A|\geq 4$ and
$L:V(\ell(G_1))\To 2^A$ an assignment for vertices of $\ell(G_1)$
such that $|L(e_5)|=2$, $|L(e_6)|=1$, $|L(e_1)|=3$, and $|L(e_i)|=4$
for $i=2,3,4$. Then $\ell(G_1)$ is $(A,L)$-colorable. (see Figure
\ref{G1}).
\end{lemma}

{\it Proof.} Let $f\in F(\ell(G_1), A)$. Suppose that
$L=L(e_4)=\{a,b,c,d\}$, $L(e_2)=L'$, $L(e_6)=\{h\}$, $f(e_4e_2)=f_1$
and $f(e_6e_2)=f_2$. If $L'\neq L-f_1$, say $a-f_1\notin L'$, we
color $e_4$ with $a$. Then color $e_6$, $e_1$, $e_5$, $e_3$ and
$e_2$ successively. So we may assume that $L'=L-f_1$. If $h-f_2\in
L'$, say $h-f_2=a-f_1$, then first color $e_4$ with $a$ and then
color $e_6$, $e_1$, $e_5$, $e_3$ and $e_2$. So $h-f_2\notin L'$ and
we color $e_6$, $e_4$,
$e_1$, $e_5$, $e_3$, $e_2$ successively.$\hfill \dashv $ \\
\begin{theorem}\label{th4-cycle}
If $G$ is a planar graph without $4$-cycles with $\D(G)\geq 5$,
then $G$ is edge-$(\D(G)+1)$-group choosable.
\end{theorem}
{\it Proof.} Let $G$ be a minimal counterexample to Theorem
\ref{th4-cycle} for some Abelian group $A$ with $|A|\geq \D(G)+1$, a
$(\D(G)+1)$-list assignment $L:V(\ell(G))\To {A\choose \D(G)+1}$ and
$f\in F(\ell(G),A)$. By Lemma \ref{lem1} and Lemma \ref{dless4}, we
may assume $\de(G)\geq 3$. Using Lemma \ref{4-cycles1}, suppose
first that $G$ contains an edge $xy$ with $d(x)+d(y)\leq 7$. Lemmas
\ref{mohem} and \ref{dless4} yield the desired contradiction. Thus
$G$ has $G_1$ as a subgraph. Remove the vertices of $\ell(G_1)$ from
$\ell(G)$, and color the remaining vertices of $\ell(G)$ from their
lists, which is possible by the minimality of $G$ as a
counterexample and using Theorem \ref{dless4}. The numbers of colors
available for vertices of $\ell(G_1)$ are at least as stated in
Lemma \ref{4-cycles2}
and so $\ell(G)$ is $(A,L,f)$-colorable, which is impossible.$\hfill \dashv $ \\

Theorem \ref{th4-cycle} together with Theorem \ref{dless4} yields
the following corollary.
\begin{corollary}
If $G$ is a planar graph without $4$-cycles, then
$\chi'_{gl}(G)\leq \D(G)+2$.
\end{corollary}

Let $G$ be a planar graph. For $v\in V$ and $f\in F$, let $m_v(f)$
denote the number of times passing through $v$ of $f$ in clockwise
order. Let $H_2$ be the subgraph induced by the edges incident with
the  $2$-vertices of $G$. It is shown that $H_2$ contains a matching
$M$ such that all $2$ vertices in $H_2$ are saturated \cite{10of5}.
If $uv\in M$ and $d(u)=2$, then $v$ is called the $2$-{\it {master}}
of $u$.
\begin{lemma}\label{charg}
If $G$ is a planar graph with maximum degree $\D(G)=4$ such that $G$
has no cycles of length from $4$ to $14$, then $G$ contains one of
the following configurations.
\begin{itemize}
\item[$(1)$] An edge $uv$ with $\min\{d(u),d(v)\}\leq 2$ and $d(u)+d(v)\leq 5$, \item[$(2)$]
A $4$-vertex $v$ incident with two $3$-faces $vv_1v_2$ and $vv_3v_4$
with $d(v_1)=d(v_4)=2$.
\end{itemize}
\end{lemma}
{\it Proof.} Let $G$ does not contain each of the mentioned
configurations. We define the initial charge function $w(x)=d(x)-4$
for each $x\in V(G)\cup F(G)$. It follows from Euler's formula that
$\sum_{x\in V(G)\cup F(G)}w(x)<0$. We define the new charge function
$\hat{w}(x)$ on $G$ as follows:
\begin{itemize}
\item[$(R_1)$] Each $r(\geq 15)$-face $f$ gives $(1-4/r)m_v(f)$ to its incident
vertex $v$ if $v$ is a cut vertex and gives 1-4/r, otherwise.
\item[$(R_2)$] Each
$2$-vertex receives $19/24$ from its neighbors if it is incident
with a $3$-face and receives $8/15$ from its $2$-master,
otherwise.\item[$(R_3)$] Each $3$-face receives $1/3$ from its
incident vertices.
\end{itemize}
It is obvious that $\hat{w}(f)=0$ for any face $f$. Let $v$ be an
arbitrary vertex of $G$. First consider the case of $d(v)=2$. If it
is incident with a $3$-face, then its other incident face $f$ must
have degree at least $16$. Since $G$ does not have configuration of
(1) any neighbor of $v$ should be of degree $4$. Hence, they can not
be $2$-vertices. It follows that $v$ receives at least $1-4/16=3/4$
from $f$ and $19/12$ from its neighbors, and gives $1/3$ to its
incident $3$-face. Otherwise, $v$ receives at least $22/15$ from its
incident faces and $8/15$ from its $2$-master. Hence,
$\hat{w}(x)\geq w(x)+\min\{3/4+19/12-1/3,22/15+8/15\}=0$. Now
consider the case of $d(v)=3$. $v$ receives at least $22/15$ from
its incident faces. Hence, $\hat{w}(v)\geq w(v)+22/15-1/3=2/15>0$.
If $d(v)=4$ and it is incident with two $3$-faces, then $v$ is
adjacent to at most one $2$-vertex since $G$ does not have
configuration of (2). It follows that $\hat{w}(v)\geq
w(v)+22/15-(2/3+19/24)=1/120>0$. Otherwise, it receives at least
$3\times 11/15$ from its incident faces, and gives at most $1/3$ to
its incident $3$-face and $19/24+8/15$ to its adjacent $2$-vertices.
It follows that $\hat{w}(v)\geq
w(v)+33/15-(1/3+19/24+8/15)=13/24>0$. This implies that $\sum_{x\in
V\cup F}w(x)=\sum_{x\in V\cup F}\hat{w}(x)>0$, a contradiction.$\hfill \dashv $ \\

Let $G$ be a planar graph with maximum degree $\D(G)=4$ such that
$G$ has no cycle of length $i$ for  $4\leq i\leq 14$. If $G$ is not
edge-$\D(G)$-group choosable, then an argument similar to the proof
of Theorem \ref{th4-cycle}, shows that $G$ does not contain any of
configurations mentioned in Lemma \ref{charg}, this leads to a
contradiction. So we have the following.
\begin{theorem}
If $G$ is a planar graph with maximum degree $\D(G)=4$ such that $G$
has no cycles of length from $4$ to $14$, then
$\chi'_{gl}(G)=\D(G)$.
\end{theorem}

In the sequel, we shall study edge-group choosability of plane
graphs without $5$-cycles, plane graphs without $5$-cycles with a
chord and plane graphs without induced $5$-cycles. A cycle $C$ of
length $k$ of a graph $G$ is called a {\it{$k$-hole}} (resp.
{\it{$k$-net}}) if $C$ has no (resp. at least one) chord in $G$. The
following is a structural lemma for plane graphs without $5$-cycles.
\begin{lemma}\cite{2}\label{5-cycles1}
If a plane graph $G$ with $\de(G)\geq 3$ has no five cycles, then
there exists an edge $xy$ of $G$ such that $d(x)=3$ and $d(y)\leq
5$.
\end{lemma}
A subgraph $H$ of a plane graph $G$ is called a {\it{cluster}} if
$H$ consists of a non-empty minimal set of $3$-faces in $G$ such
that no other $3$-face is adjacent to a member of this set. A
structural lemma for plane graphs with $\D(G)=6$ and without
$5$-nets is as follows.
\begin{lemma}\cite{7}\label{5-cycles2}
Let $G$ be a planar graph with $\D(G)=6$ and without $5$-nets. Then
$G$ contains one of the following configurations.
\begin{itemize}
\item[$(1)$] An edge $xy$ with $d(x)+d(y)\leq 8$, \item[$(2)$] A
$4$-cycle $C=uvwx$ such that $d(u)=d(w)=3$ and $d(v)=d(x)=6$,
\item[$(3)$] The subgraph $G_2$ consisting of a $6$-vertex $x$
incident to four  $3$-faces $[xx_1x_2]$, $[xx_2x_3]$, $[xx_4x_5]$,
$[xx_5x_6]$ such that $d(x_1)=d(x_5)=3$, $d(x_3)=4$ and
$d(x_2)=d(x_4)=d(x_6)=6$ (see Figure 2), \item[$(4)$] The graph
$G_3$ consisting of a cluster with $\{u,x,y,z\}$ as its vertex set
such that $d(u)=3$, $d(x)=d(y)=d(z)=6$, and satisfying the following
properties: \\$(a)$ $x$ is incident to $(3,6,6)$-faces $[uxy]$,
$[uxz]$ and $[xx_1x_2]$ with $d(x_1)=3$, \\$(b)$ $y$ is incident to
$(3,6,6)$-faces $[uxy]$, $[yuz]$, $[yy_1y_2]$ and $[yy_2y_3]$ with
$d(y_2)=3$.
\end{itemize}
\end{lemma}

For planar graphs without $5$-nets, we need the following lemma.
\begin{lemma} \cite{7}\label{5-cycles3}
Every planar graph $G$ with $\de(G)\geq 3$ and without $5$-nets
contains an edge $xy$ such that $d(x)+d(y)\leq 9$.
\end{lemma}

For studying Conjecture $3$ for planar graphs with maximum degree
$5$ and without $5$-nets we need the following structural lemma.
\begin{lemma} \cite{7} \label{5-cycles4,5}
Consider the following configurations.
\begin{itemize}
\item[$(1)$] An edge $xy$ with $d(x)+ d(y)\leq 7$,\item[$(2)$] A
$4$-cycle $vwxu$ such that $d(u)=d(w)=3$ and $d(v)=d(x)=5$,\item[
$(3)$] The subgraph $G_4$ consisting of a $5$-vertex incident to
three $(3,5,5)$-faces,\item[$(4)$] The subgraph $G_5$ consisting of
a $4$-vertex $u$ adjacent to three $4$-vertices $x, y, z$ and
incident to a $3$-face $[uxy]$.
\end{itemize}
Let $G$ be a planar graph with $\D(G)=5$. Then \\$(a)$ If $G$ is
without $5$-nets and $6$-nets, then it contains at least one of the
graphs mentioned in $(1)$, $(2)$ and $(3)$, \\$(b)$ If $G$ is
without $4$-nets and $5$-nets, then it contains at least one of the
graphs mentioned in $(1)$, $(2)$ and $(4)$.
\end{lemma}
\begin{lemma}\cite{10}\label{5-cycles6}
Let $G$ be a connected planar graph with $\delta(G)\geq 2$. If $G$
contains no $5$-cycles or no $6$-cycles, then $G$ contains either
an edge $xy$ with $d(x)+d(y)\leq 9$ or a $2$-alternating cycle.
\end{lemma}
\begin{theorem}\label{th5-cycles1}
If $G$ is a planar graph without $5$-cycles with maximum degree
$\D$, then $G$ is edge-$(\D+2)$-group choosable.
\end{theorem}
{\it Proof.} Let $G$ be a minimal counterexample to Theorem
\ref{th5-cycles1}. Then by Lemma \ref{lem1} and Theorem
\ref{dless4} we have $\de(G)\geq 3$,  $\D\geq 5$ and by Lemma
\ref{5-cycles1}, there exists a vertex $e\in V(\ell(G))$ with
$d_{\ell(G)}(e)\leq 6$, which is impossible by Lemma \ref{mohem}.$\hfill \dashv $ \\

\begin{figure}[t]
  \begin{center}
  \includegraphics[width=4cm]{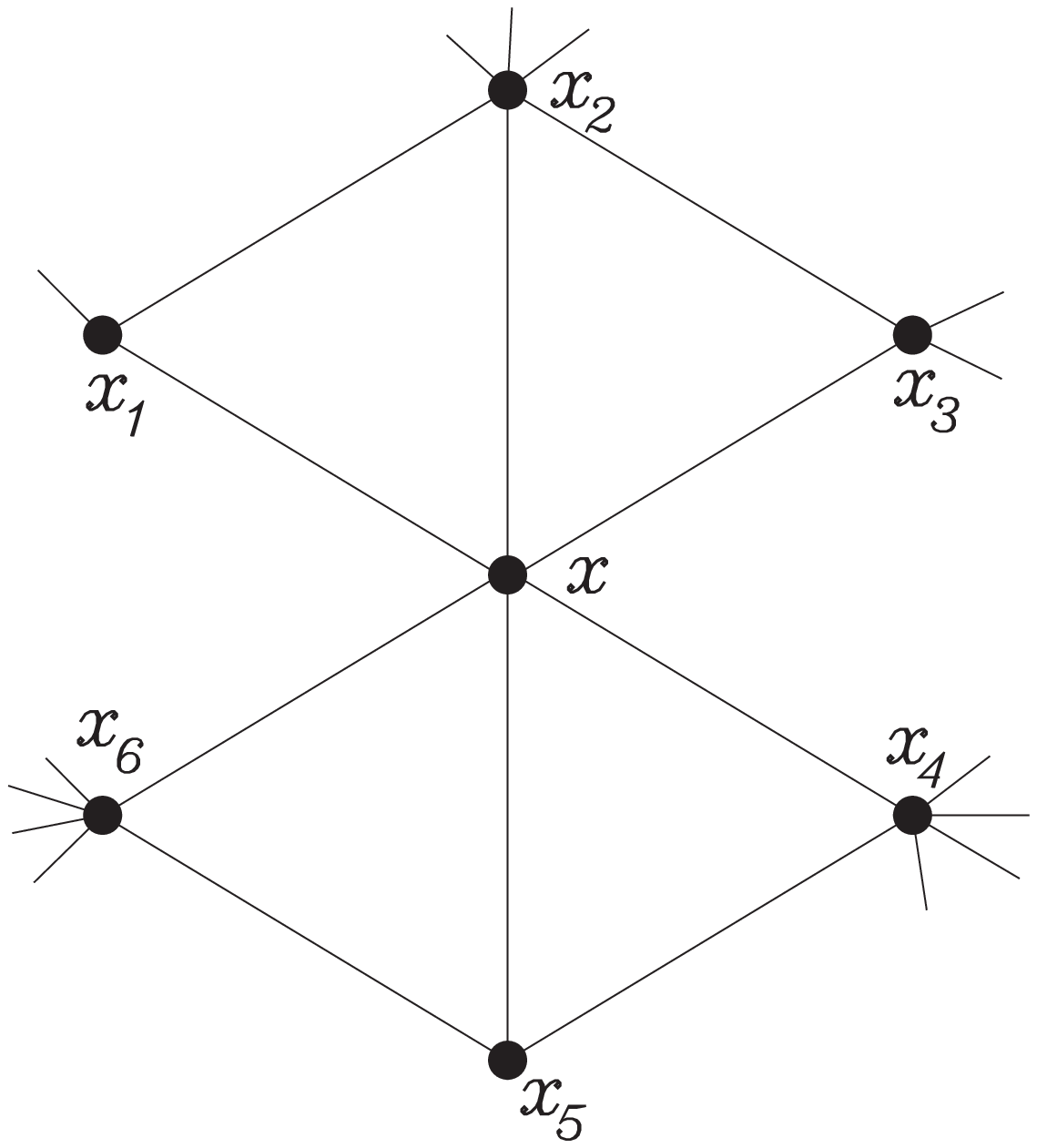}\hspace*{2cm}
  \includegraphics[width=3cm]{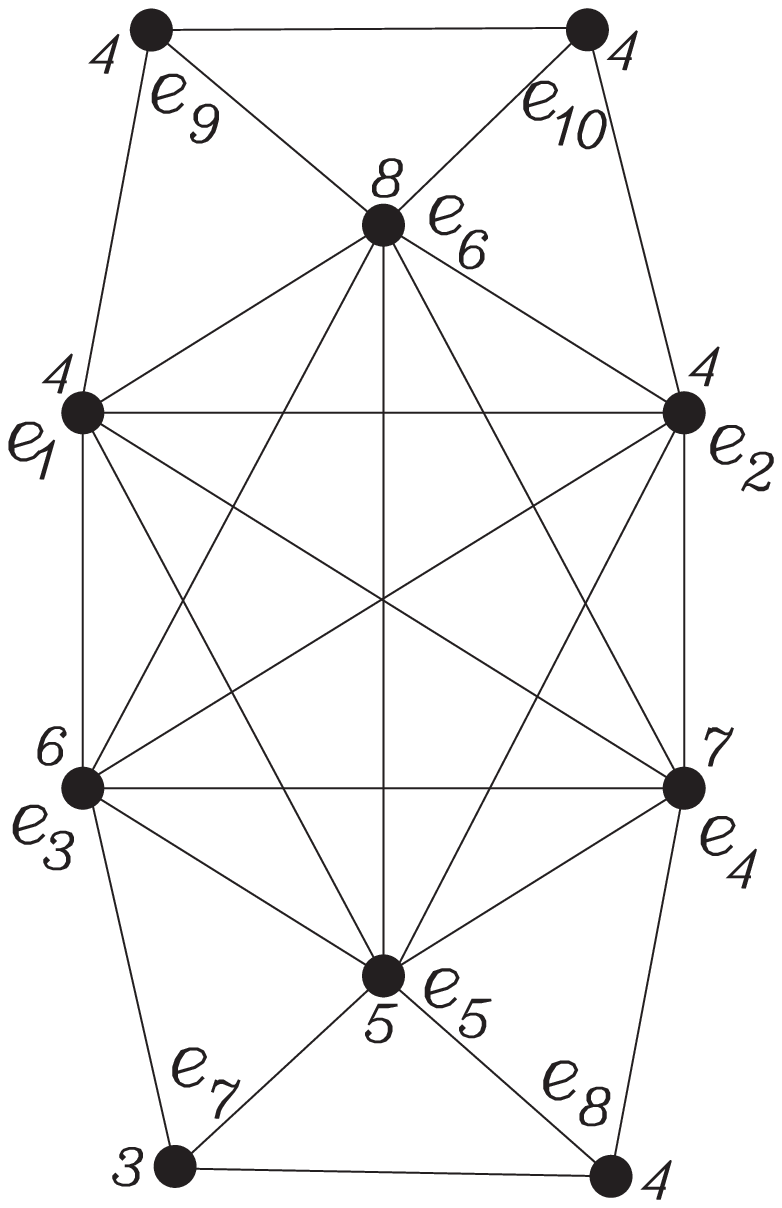}\\
  \caption{The graph $G_2$ and its line graph.}\label{G2}
  \end{center}
\end{figure}

For graphs with maximum degree at least $7$, we have a stronger
result as follows.
\begin{theorem}\label{th5-cycles2}
If $G$ is a planar graph with maximum degree $\D\geq 7$ and
without $5$-cycles, then $\chi'_{gl}(G)\leq \D+1$.
\end{theorem}
{\it Proof.} Let $G$ be a minimal counterexample to Theorem
\ref{th5-cycles2} for some Abelian group $A$ with $|A|\geq \D(G)+1$,
a $(\D(G)+1)$-list assignment $L:V(\ell(G))\To {A\choose \D(G)+1}$
and $f\in F(\ell(G),A)$. Then by Theorem \ref{th5-cycles1} and Lemma
\ref{lem1}, $\de(G)\geq 3$ and we can apply Lemma \ref{5-cycles6}.
Hence there exists a vertex $e\in V(\ell(G))$ with
$d_{\ell(G)}(e)\leq 7$, which is impossible by
Theorem \ref{th5-cycles1} and Lemma \ref{mohem}. $\hfill \dashv $ \\
\begin{theorem}\label{th5-cycles3} Let $G$ be a planar graph without
$5$-nets,
\begin{itemize}
\item[$(1)$] If $\D(G)\geq 6$, then $\chi'_{gl}(G)\leq
\max\{8,\D(G)+1\}$, \item[$(2)$] If $\D(G)=5$ and $G$ contains no
$4$-nets or $6$-nets, then $\chi'_{gl}(G)\leq 7$.
\end{itemize}
\end{theorem}
{\it Proof.} $(1)$ Let $k=\max\{8, \D(G)+1\}$ and $G$ be a minimal
counterexample to Theorem \ref{th5-cycles3}. Then there is an
Abelian group $A$ with $|A|\geq k$, a $k$-assignment
$L:V(\ell(G))\To {A\choose k}$ and $f\in F(\ell(G), A)$ such that
$\ell(G)$ is not $(A, L, f)$-colorable. Suppose first that
$\D(G)=6$. By Lemma \ref{5-cycles2}, we consider four cases as
follows. \\Cases $1$, $2$. $G$ contains an edge $xy$ with
$d(x)+d(y)\leq 8$ or a $4$-cycle $C=uvwx$ such that $d(u)=d(w)=3$
and $d(v)=d(x)=6$. Both cases lead to a contradiction by an argument
similar to the proofs of Lemma \ref{lem1} and Theorem
\ref{th5-cycles2}. \\Case $3$. $G$ contains graph $G_2$ of Lemma
\ref{5-cycles2} (see Figure \ref{G2}).
Remove the vertices of $\ell(G_2)$ from $\ell(G)$ and color the
remaining vertices of $\ell(G)$ from their lists, which is possible
by the minimality of $G$ as a counterexample. There are now $4, 4,
6, 7, 5, 8, 3, 4, 4, 4$ colors available for the vertices
$e_1,\ldots,e_{10}$, respectively where $e_1=xx_4$, $e_2=xx_6$,
$e_3=xx_3$, $e_4=xx_1$, $e_5=xx_2$, $e_6=xx_5$, $e_7=x_2x_3$,
$e_8=x_1x_2$, $e_9=x_4x_5$, and $e_{10}=x_5x_6$. It is easily seen
that we can color these vertices in their order. Thus $G$ is not a
counterexample, which is a contradiction. \\Case $4$. $G$ contains
graph $G_3$ of Lemma \ref{5-cycles2}. The
proof is similar to case $3$.\\
If $\D(G)>6$, by applying Lemmas \ref{lem1} and  \ref{5-cycles3} the
assertion holds.

$(2)$ If $\D(G)=5$ and $G$ has no $4$-nets, apply Lemma
\ref{5-cycles4,5}$(b)$ and
if $G$ has no $6$-nets, apply Lemma \ref{5-cycles4,5}$(a)$.$\hfill \dashv $ \\

The structure of planar graphs without non-induced $5$-cycles is
given in the following lemma.
\begin{lemma} \cite{8}\label{5-cycles7}
Let $G$ be a planar graph without non-induced $5$-cycles. Then $G$
contains one of the following configurations.
\begin{itemize}
\item[$(1)$] An edge $uv$ with $d(u)+ d(v)\leq
\max\{8,\D(G)+2\}$,\item[$(2)$] An even cycle $C: v_1v_2 \ldots
v_{2n}$ with $d(v_1) = d(v_3) = \cdots = d(v_{2n-1}) = 3$ and
$d(v_2) = d(v_4) = \cdots = d(v_{2n}) = \D(G)$.
\end{itemize}
\end{lemma}
\begin{theorem}\label{th5-cycles4}
If $G$ is a planar graph without non-induced $5$-cycles, then
$\chi'_{gl}(G)\leq \max\{7, \D+2\}$.
\end{theorem}
{\it Proof.} Applying Lemma \ref{5-cycles7}, the proof is similar
to the proof of Theorem \ref{th5-cycles2}.$\hfill \dashv $ \\

If planar graph $G$ without non-induced $5$-cycles in addition
contains no even cycles, we can replace $\D+2$ by $\D+1$ in Theorem
\ref{th5-cycles4}.
\begin{corollary}\label{th5-cycles5}
If $G$ is a planar graph without non-induced $5$-cycles and without
even cycle $C: v_1v_2 \ldots v_{2n}$ with $d(v_1) = d(v_3) = \cdots
= d(v_{2n-1}) = 3$ and $d(v_2) = d(v_4) = \cdots = d(v_{2n}) =
\D(G)$, then $\chi'_{gl}(G)\leq \max\{7, \D+1\}$.
\end{corollary}

We finish this section by studying Conjecture $3$ on planar graphs
without $3$-, $6$- or $7$-cycles. Because the proofs are similar to
the proofs we did before, we omit them. The structural lemma for
triangle free planar  graphs is as follows.
\begin{lemma}\cite{3}\label{3-cycles1}
Let $G$ be a triangle-free plane graph with $\D(G)\geq 5$. If
$d(x)+d(y)\geq \D(G)+3$ for every edge $xy\in E(G)$, then $\D(G)=
5$ and $G$ contains a $4$-face incident with two $3$-vertices and
two $5$-vertices.
\end{lemma}

The following lemmas show the structure of planar graphs without
adjacent triangles and without $6$-cycles, respectively.
\begin{lemma}\cite{9}\label{3-cycles2}
Let $G$ be a plane graph without adjacent triangles. Then $G$
contains one of the following configurations.
\begin{itemize}
\item[$(1)$] An edge $uv$ with $d(u)+ d(v)\leq \max\{8, \D+2\}$,
\item[$(2)$] A $4$-cycle $vwxu$ such that $d(u)=d(w)=3$ and
$d(v)=d(x)=\D(G)$, \item[$(3)$] The subgraph $G_6$ consisting of a
$6$-vertex $v$ incident with two $(3, 6, 6)$-faces and one $(4, 5,
6)$-face.
\end{itemize}
\end{lemma}
\begin{lemma}\cite{11}\label{6-cycles1}
If $G$ is a plane graph without $6$-cycles and $\de(G) \geq 3$,
then there is an edge $xy \in E(G)$ such that $d(x) + d(y) \leq
8$.
\end{lemma}
\begin{lemma}\cite{9}\label{7-cycles1}
Every planar graph $G$ without $7$-cycles contains one of the
following configurations.
\begin{itemize}
\item[$(1)$] An edge $uv$ with $d(u)+d(v)\leq  \max\{9,
\D(G)+2\}$, \item[$(2)$] A $4$-cycle $vwxu$ such that
$d(u)=d(w)=3$ and $d(v)=d(x)=\D(G)$.
\end{itemize}
\end{lemma}

Using Lemmas \ref{lem1}, \ref{3-cycles1}-\ref{7-cycles1} and
arguments similar to the proofs we did before, we can summarize all
results on Conjecture $3$ for planar graphs without $3$-, $6$- or
$7$-cycles in the following theorem.
\begin{theorem}\label{th367-cycles1}
Let $G$ be a planar graph,
\begin{itemize}
\item[$(1)$] If $G$ contains no $3$-cycles and $\D(G)\geq 6$, then
$\chi'_{gl}(G)\leq \D(G)+1$, \item[$(2)$] If $G$ contains no
$3$-cycles and $\D(G)=5$, then $\chi'_{gl}(G)\leq 7$, \item[$(3)$]
If $G$ contains no adjacent triangles, then $\chi'_{gl}(G)\leq
\max\{7, \D(G)+2\}$, \item[$(4)$] If $G$ contains no $6$-cycles,
then $\chi'_{gl}(G)\leq \max\{7, \D(G)+1\}$, \item[$(5)$] If $G$
contains no $7$-cycles, then $\chi'_{gl}(G)\leq \max\{8, \D(G)+2\}$.
\end{itemize}
\end{theorem}

As an immediate consequence of Theorem \ref{th367-cycles1} and
Theorem \ref{dless4}, we have the following.
\begin{corollary}\label{th6-cycles2}
If $G$ is a planar graph without $6$-cycles or without $3$-cycles,
then $\chi'_{gl}(G)\leq  \D(G)+2$.
\end{corollary}
\section{\normalsize\bf{Edge-group choosability of outerplanar and near-outerplanar
graphs}}\label{secouterplanar} A planar graph is called
{\it{outerplanar}} if it has a drawing in which each vertex lies on
the boundary of the outer face. It is well-known that a graph is
outerplanar iff it contains neither $K_4$ nor $K_{2, 3}$ as a minor
(see for example \cite{diestel}).
In this section, we see that Conjecture $3$ is true for outerplanar
graphs and {\it{near-outerplanar}} graphs, i.e. graphs that are
either $K_4$-minor- free or $K_{2, 3}$-minor-free.
\begin{lemma}\cite{Thirteen}\label{outerplanar1}
If $G$ is an outerplanar graph, then at least one of the following
cases holds.
\begin{itemize}
\item[$(1)$] $\de(G)=1$, \item[$(2)$] There exists an edge $uv$
such that $d(u)=d(v)=2$, \item[$(3)$] There exists a $3$-face $uxy$
such that $d(u)=2$, $d(x)=3$, \item[$(4)$] $G$ contains the graph
$G_7$ consisting of two $3$-faces $xu_1v_1$ and $xu_2v_2$ such that
$d(u_1)=d(u_2)=2$ and $d(x)=4$ and these five vertices are all
distinct.
\end{itemize}
\end{lemma}

This lemma lets us to prove that Conjecture $3$ is true for
outerplanar graphs. In fact, we prove a stronger result.
\begin{theorem}\label{thouterplanar1}
Let $G$ be an outerplanar graph with maximum degree $\D$. Then
\begin{itemize}
\item[$(1)$] If $\D<5$, then $\chi'_{gl}(G)\leq \D+1$,
\item[$(2)$] If $\D\geq 5$, then $\chi'_{gl}(G)\leq \D$.
\end{itemize}
\end{theorem}
{\it Proof.} We prove the second part of the theorem, the proof of
the first part is similar. Let $G$ be a minimal counterexample to
the second part of Theorem \ref{thouterplanar1}. Then there is an
Abelian group $A$ with $|A|\geq \D$, a $\D$-assignment
$L:V(\ell(G))\To {A\choose \D}$ and $f\in F(\ell(G), A)$ such that
$\ell(G)$ is not $(A, L, f)$-colorable. By Lemma \ref{outerplanar1},
we consider four cases as follows.

Case $1$. $\de(G)=1$. Then $G$ contains an edge $uv$ such that
$d(u)=1$ and $d(v)\leq \D$. Hence $\ell(G)$ contains a vertex $e=uv$
with $d_{\ell(G)}(e) \leq \D-1$, a contradiction by the first part
of this theorem and Lemma \ref{mohem}.

Case $2$. There exists an edge $uv$ such that $d(u)=d(v)=2$. Then
$\ell(G)$ contains a vertex $e$ with $d_{\ell(G)}(e)=2$, a
contradiction by the first part of this theorem and Lemma
\ref{mohem}.

Case $3$. There exists a $3$-face $C : uxy$ such that $d(u)=2$,
$d(x)=3$. Remove the vertices of $\ell(C)$ from $\ell(G)$ and color
the remaining vertices of $\ell(G)$ from their lists. There are now
at least $1, 2, 4$ colors available for the vertices of $\ell(C)$.
It is easily seen that we can color these vertices. Thus $G$ is not
a counterexample which is a contradiction.

\begin{figure}[h]
  \begin{center}
  \includegraphics[width=3cm]{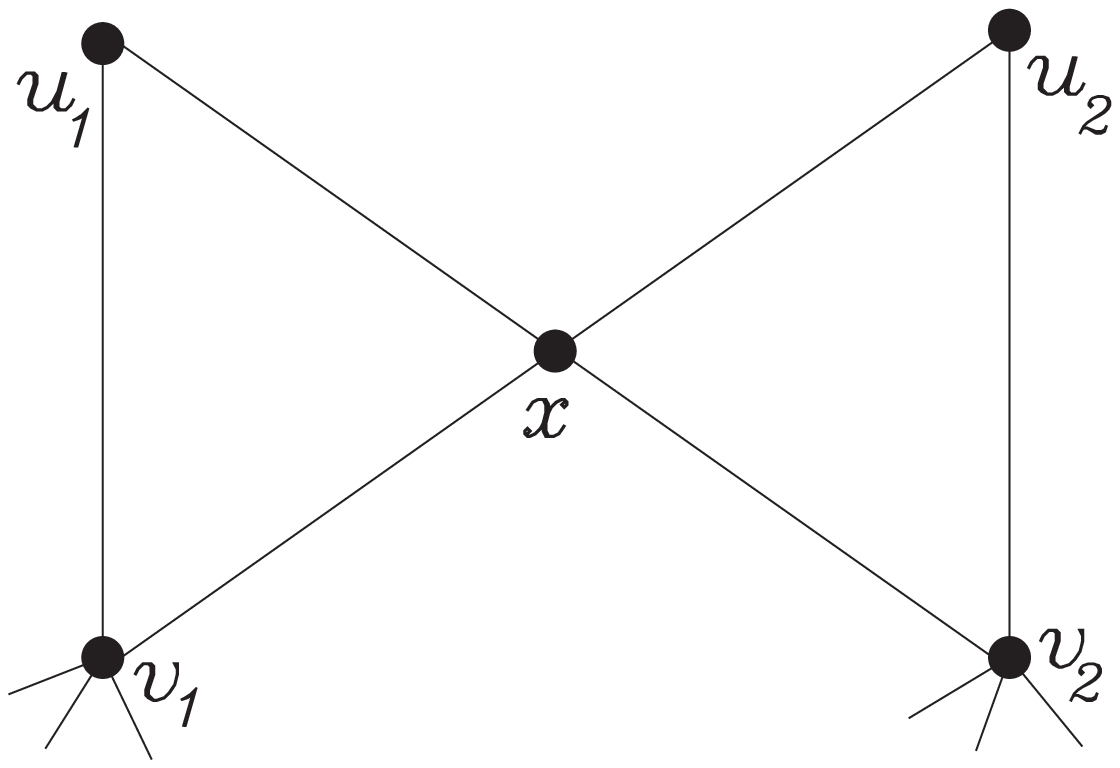}\hspace*{2cm}
  \includegraphics[width=4.5cm]{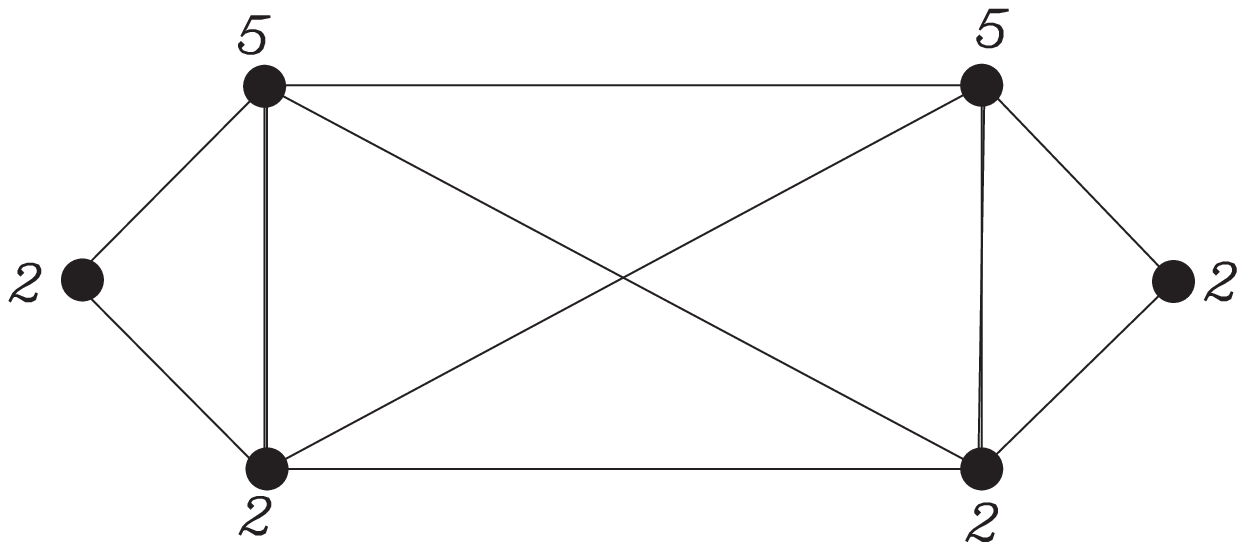}\\
  \caption{The graph $G_7$  and its line graph.}\label{G7}
  \end{center}
\end{figure}

Case $4$. $G$ contains the subgraph $G_7$ of Lemma
\ref{outerplanar1} (see Figure \ref{G7}).
Remove the vertices of $\ell(G_7)$ from $\ell(G)$ and color the
remaining vertices of $\ell(G)$ from their lists, which is possible
by the minimality of $G$ as a counterexample and using the second
part of this theorem. There are now at least $ 2, 2, 2, 2, 5, 5$
colors available for the vertices of $\ell(G_7)$. It is easily seen
that we can color these vertices in the above order. Thus $G$ is not
a counterexample which is a
contradiction.$\hfill \dashv $ \\

A graph is called {\it{series-parallel}} if it has no subgraph
isomorphic to a subdivision of $K_4$. It is well-known
\cite{4of12} that every simple series-parallel graph has a vertex
of degree at most two. Thus using Lemma \ref{lem1}, we have the
following corollary.
\begin{corollary}\label{thserie-parallel}
If $G$ is a simple series-parallel graph, then $\chi'_{gl}(G)\leq
\D(G)+1$. In particular, every $K_4$-minor-free graph $G$ is
edge-$(\D(G)+1)$-group choosable.
\end{corollary}

In the following, we will prove that Conjecture $3$ holds for
$({K_2}^c+(K_1 \cup K_2))$-minor-free graphs, where ${K_2}^c+(K_1
\cup K_2)$ is the graph obtained from the union of $K_1$ and $K_2$
and joining them by ${K_2}^c$, or, equivalently, it is the graph
obtained from $K_{2,3}$ by adding an edge joining two vertices of
degree $2$. This implies that Conjecture $3$ is true for
$K_{2,3}$-minor-free graphs.
\begin{lemma}\cite{13}\label{near-outerplanar1}
Let $G$ be a $({K_2}^c+(K_1 \cup K_2))$-minor-free graph. Then each
block of $G$ is either $K_4$-minor-free or isomorphic to $K_4$.
\end{lemma}

A graph $G$ is called {\it{$D$-group choosable}} if it is $(A,
L)$-colorable for every Abelian group $A$ with $|A|\geq \D(G)$ and
every list assignment $L : V (G) \To 2^A$ with $|L(v)|=d(v)$ for
each vertex $v$. There is a characterization of all $D$-group
choosable graphs in \cite{main} as follows.
\begin{theorem}\label{D-group}
A connected graph $G$ is not $D$-group choosable iff every block
of $G$ is either complete or cycle.
\end{theorem}

By Theorem \ref{D-group},  it is easily seen that $\ell(K_4)$ is
$D$-group choosable and so $\chi'_{gl}(K_4)\leq 4$. Moreover, it is
well-known \cite{10of13} that a $K_4$-minor-free graph $G$ with $|V
(G)|\geq 4$ has at least two non-adjacent vertices with degree at
most $2$. By modifing some proofs in \cite{13}, we have the
following result.
\begin{theorem}\label{thnear-outerplanar1}
If $G$ is a $({K_2}^c+(K_1 \cup K_2))$-minor-free graph, then
$\chi'_{gl}(G)\leq \D(G)+1$. In particular, every $K_{2,
3}$-minor-free graph $G$ is edge-$(\D(G)+1)$-group choosable.
\end{theorem}
{\it Proof.} Let $G$ be a minimal counterexample to Theorem
\ref{thnear-outerplanar1}. Then there is an Abelian group $A$ with
$|A|\geq \D(G)+1$, a $(\D(G)+1)$-assignment $L:V(\ell(G))\To
{A\choose \D(G)+1}$ and $f\in F(\ell(G), A)$ such that $\ell(G)$ is
not $(A, L, f)$-colorable. Clearly $G$ is connected, $G\neq K_4$ and
by Lemma \ref{lem1}, $\de(G)\geq 3$. If $\D(G)=3$, then $G$ is
$3$-regular and not complete, a contradiction follows since
$\ell(G)$ is $4$-regular and so it is $D$-group choosable by Theorem
\ref{D-group}. So we may assume that $\D(G)\geq 4$. By Lemma
\ref{near-outerplanar1}, every block of $G$ is either
$K_4$-minor-free or isomorphic to $K_4$. We first show that $G$ is
not $2$-connected. Suppose on the contrary that $G$ is
$2$-connected. Then $G$ is $K_4$-minor free since $\D(G)\geq 4$ and
so by Corollary \ref{thserie-parallel}, $\chi'_{gl}(G)\leq \D(G)+1$,
a contradiction. Thus $G$ is not $2$-connected. Let $B$ be an
end-block of $G$ with cut-vertex $z_0$. Clearly $B\ncong K_2$.
Moreover, if $B\cong K_4$, it is easily seen that an $(A, L,
f)$-coloring of $\ell(G)-V(\ell(B))$ can be extended to an $(A, L,
f)$-coloring of $\ell(G)$. This contradiction shows that $B\ncong
K_4$. Hence using Lemma \ref{near-outerplanar1}, $B$ is
$K_4$-minor-free, and so it contains at least two vertices of degree
at most $2$. Hence $G$ has at least one vertex of degree at most
$2$, this contradiction completes the proof of the theorem.$\hfill
\dashv $ 
\vspace{1.5cm}

 \end{document}